\documentclass{elsart}

\usepackage{bm}
\usepackage{amsmath}
\usepackage{amsfonts}
\usepackage{amssymb}

\newtheorem{theorem}{Theorem}[section]

\newtheorem{proposition}[theorem]{Proposition}
\newtheorem{corollary}[theorem]{Corollary}
\theoremstyle{definition}

\newtheorem{remark}[theorem]{Remark}
\newtheorem{example}[theorem]{Example}

\newcommand{\RR}{{\mathbb R}}

\def \bv{\mathbf{v}}
\def \x{\mathbf{x}}
\def \n{\mathbf{n}}

\def \z{\mathbf{z}}

\def \e{\mathbf{e}}

\def \R{\mathbf{R}}

\def \H{\mathbf{H}}

\def \Q{\mathbf{Q}}

\def \I{\mathbf{I}}

\def \E{\mathbf{E}}

\def \S{\mathbf{S}}
\def \P{\mathbf{P}}

\def \L{\mathbf{L}}
\def \F{\mathbf{F}}

\def \0{\mathbf{0}}
\def \tr{\textnormal{tr}}

\begin{document}

\begin{frontmatter}

\title{Robust Dimension Reduction, Fusion Frames, and Grassmannian Packings}

\author{Gitta Kutyniok\corauthref{cor1}\thanksref{GK}},
\ead{kutyniok@math.princeton.edu}
\author{Ali Pezeshki\thanksref{AP}},
\ead{pezeshki@math.princeton.edu}
\author{Robert Calderbank\thanksref{AP}},
\ead{calderbk@math.princeton.edu}
\author{Taotao Liu\thanksref{TL}}
\ead{taotaol@math.princeton.edu}

\address{Program in Applied and Computational Mathematics, Princeton University, Princeton, NJ 08544-1000, USA}

\corauth[cor1]{Corresponding author}
\thanks[GK]{Supported by Deutsche Forschungsgemeinschaft (DFG)
Heisenberg-Fellowship KU 1446/8-1.}
\thanks[AP]{Supported in part by AFOSR Grant 00852833.}
\thanks[TL]{Supported by PACM Higgins Research Fund.}

\begin{abstract}
We consider estimating a random vector from its noisy projections
onto low-di\-men\-sional subspaces constituting a \textit{fusion
frame}. A fusion frame is a collection of subspaces, for which the
sum of the projection operators onto the subspaces is bounded
below and above by constant multiples of the identity operator. We
first determine the minimum mean-squared error (MSE) in linearly
estimating the random vector of interest from its fusion frame
projections, in the presence of white noise. We show that MSE
assumes its minimum value when the fusion frame is \textit{tight}.
We then analyze the robustness of the constructed linear minimum
MSE (LMMSE) estimator to erasures of the fusion frame subspaces.
We prove that tight fusion frames consisting of equi-dimensional
subspaces have maximum robustness (in the MSE sense) with respect
to erasures of one subspace, and that the optimal subspace
dimension depends on signal-to-noise ratio (SNR). We also prove
that tight fusion frames consisting of equi-dimensional subspaces
with equal pairwise chordal distances are most robust with respect
to two and more subspace erasures. We call such fusion frames {\em
equi-distance tight fusion frames,} and prove that the chordal
distance between subspaces in such fusion frames meets the
so-called {\em simplex bound,} and thereby establish connections
between equi-distance tight fusion frames and {\em optimal
Grassmannian packings.} Finally, we present several examples for
construction of equi-distance tight fusion frames.
\end{abstract}


\end{frontmatter}

\section{Introduction}

The notion of a {\em fusion frame} (or {\em frame of subspaces})
was introduced by Casazza and Kutyniok in \cite{CK04} and further
developed by Casazza et al. in \cite{CKL07}. A fusion frame for
$\mathbb{R}^{M}$ is a finite collection of subspaces
$\{\mathcal{W}_i\}_{i=1}^{N}$ in $\mathbb{R}^{M}$ such that there
exist constants $0 < A \le B < \infty$ satisfying
\[
A\|\x\|^2 \le \sum_{i=1}^N \|\P_i \x\|^2 \le B\|\x\|^2,\quad
\textnormal{for any $\x \in \mathbb{R}^{M}$},
\]
where $\P_i$ is the orthogonal projection onto $\mathcal{W}_i$.
Alternatively, $\{\mathcal{W}_i \}_{i=1}^N$ is a fusion frame if
and only if
\begin{equation}\label{eq:ff2}
A \, \I \le \sum_{i=1}^N \P_i \le B \, \I.
\end{equation}
The constants $A$ and $B$ are called {\em (fusion) frame bounds}.
An important class of fusion frames are {\em tight fusion frames},
for which $A$ and $B$ can be chosen to be equal and $\sum_{i=1}^N
\P_i = A \, \I$. We note that the definition given in \cite{CK04}
and \cite{CKL07} for fusion frames apply to closed and weighted
subspaces in any Hilbert space. However, since the scope of this
paper is limited to non-weighted subspaces in $\mathbb{R}^{M}$,
the definition of a fusion frame is only presented for this case.

A fusion frame can be viewed as a frame-like collection of
low-dimensional subspaces. In frame theory, an input signal is
represented by a collection of {\em scalars}, which measure the
magnitudes of the projections of the signal onto frame vectors,
whereas in fusion frame theory an input signal is represented by a
collection of {\em vectors}, which are the projections of the
signal onto the fusion frame subspaces. Similar to frames, fusion
frames can be used to provide a redundant and non-unique
representation of a signal. In fact, in many applications, where
data has to be processed in a distributed manner by combining
several locally processed data vectors, fusion frames can provide
a more natural mathematical framework than frames. A few examples
of such applications are as follows.

\textit{Distributed sensing:} In distributed sensing, typically a
large number of inexpensive sensors are deployed in an area to
measure a physical quantity such as temperature, sound, vibration,
pressure, etc., or to keep an area under surveillance for target
detection and tracking. Due to practical and economical factors,
such as low communication bandwidth, limited signal processing
power, limited battery life, or the topology of the surveillance
area, the sensors are typically deployed in clusters, where each
cluster includes a unit with higher computational and transmission
power for local data processing. Thus, a typical large sensor
network can be viewed as a redundant collection of subnetworks–
forming a set of subspaces.
%The primary goal is to have local measurements transmitted to a local
%sub-station within the subspace for subspace combining.
The gathered subspace information is submitted to a central
processing station for joint processing. Some references that
consider fusion frames for distributed sensing are
\cite{RGJ06},\cite{RJ06}, and \cite{CKLR07}.

\textit{Parallel Processing:} If a frame system is simply too
large to handle effectively (from the numerical stand point), we
can divide it into multiple small subsystems for simpler and parallel
processing. By introducing redundancy, when splitting the large
system, we can introduce robustness against errors due to failure
of a subsystem. Fusion frames provide a natural framework for
splitting a large frame system into smaller subsystems and then
recombining the subsystems.

\textit{Packet encoding:} In digital media transmission,
information bearing source symbols are typically encoded into a
number of packets and then transmitted over a communication
network, e.g., the internet. The transmitted packet may be
corrupted during the transmission or completely lost due to, for
example, buffer overflows. By introducing redundancy in encoding
the symbols, according to an error correcting scheme, we can
increase the reliability of the communication scheme. Fusion
frames, as redundant collections of subspaces, can be used to
produce a redundant representation of a source symbol. In the
simplest form, we can think of each low-dimensional projection as
a packet that carries some new information about the symbol. At
the destination the packets can be decoded jointly to recover the
transmitted symbol. The use of fusion frames for packet encoding
is considered in \cite{Bod07}.

In this paper, we consider estimating a \textit{random vector}
from its noisy projections onto low-dimensional subspaces
constituting a fusion frame. As far as we know, optimal
reconstruction of random vector signals from fusion frame
measurements (or even frame measurements) has not been considered
before, despite the fact that random vectors provide a natural way
of modeling signals in many applications, including distributed
sensing.

The optimal reconstruction of a random signal is different from
the canonical reconstruction of a deterministic signal in a fusion
frame that is considered in \cite{CKL07}. The canonical
reconstruction strategy of a deterministic signal
$\x\in\mathbb{R}^M$ from its fusion frame measurements involves
using the fusion frame operator $\S \x = \sum_{i=1}^N \P_i \x$,
which is invertible and self-adjoint. The deterministic signal
$\x$ can then be recovered from its measurements $\{\P_i \x\}_{i =
1}^N$ as $\x = \S^{-1} \S \x = \sum_{i=1}^N \S^{-1} \P_i \x$.
However, this strategy is not optimal in the MSE sense if $\x$ is
a random vector.

When $\x$ is random, the (linear) strategy that achieves minimum
MSE is {\em linear minimum mean-squared error (LMMSE) estimation} or
{\em Wiener filtering}, which is well-known in the statistical signal
processing literature (cf. \cite[Ch.8]{Scharf-91}). We use this
strategy to estimate the random vector $\x$ -- assuming that its
covariance matrix $\R_{xx}=E[\x\x^T]$ is nonsingular and known
(but otherwise arbitrary) -- from its noisy fusion frame
projections. We determine the MSE in the LMMSE estimation of $\x$
and show that the MSE assumes its minimum value when the fusion
frame is \textit{tight}. Our analysis also clarifies the effect of
additive white noise on signal estimation in fusion frames, which
has not been studied before.

We then analyze the robustness of LMMSE estimators to erasures of
the fusion frame subspaces. Erasures of subspaces can occur due to
many factors in practice. In the distributed sensing example, a
subspace erasure can occur due to a faulty or out of battery
cluster of sensors, or due to loss of data during the transmission
of local subspace information to the central processor. In
scenarios where one or more sensor clusters are believed to be out
of range for measuring the signal, or blocked by obstacles, their
corresponding subspaces can be discarded on purpose. In the
parallel processing example, an erasure can occur when a local
processor crashes. In the packet encoding example, an erasure can
occur when buffers in the network overflow.

Constructing frames that allow for robust reconstruction of a
\textit{deterministic signal} in the presence of frame element
erasures has been considered by a number of authors. In
\cite{GKK01}, Goyal et al. show that a normalized frame is
optimally robust against noise and one erasure (erasure of one
element of the frame) if the frame is tight. Some ideas
concerning multiple erasures were also presented. The work of Casazza and
Kova{\v{c}}evi{\'c} \cite{CK03} focuses mainly on designing frames,
which maintain completeness under a particular number of erasures.
Holmes and Paulsen \cite{HP04} and Bodmann and Paulsen \cite{BP05}
study the robustness of frames under multiple erasures and show that
maximal robustness with respect to the worst-case (maximum)
Euclidean reconstruction error is achieved when the frame elements
are equi-angular. The connection between equi-angular frames and
equi-angular lines has also been explored by Strohmer and Heath in
\cite{SH03}, where the so-called Grassmannian frames are introduced.

There are also a few papers that consider the construction of
fusion frames for robust reconstruction of \textit{deterministic
signals} in the presence of subspace erasures. The main result in
this context is due to Bodmann \cite{Bod07}, who shows that a tight
fusion frame is optimally robust against one subspace erasure
if the dimensions of the subspaces are equal. He also proves that
a tight fusion frame is optimally robust against multiple erasures
if the subspaces satisfy the so-called \textit{equi-isoclinic}
condition. The performance measure considered in \cite{Bod07} is
the worst-case (maximum) Euclidean reconstruction error. The
equi-isoclinic condition requires all pairs of subspaces to have
the same set of principal angles. This condition is very
restrictive and there are only a few known examples of fusion
frames that satisfy this condition. The single erasure case
discussed in \cite{Bod07} has also been studied by Casazza and
Kutyniok in \cite{CK07}. We emphasize that all the above work on
robustness to erasures in frames and fusion frames deals with the
case where the signal of interest is deterministic.

In this paper, we analyze the effect of subspace erasures on the
performance of LMMSE estimators. We determine how the MSE of an
LMMSE estimator, constructed based on the second-order statistics
of the data in the absence of erasures, is affected by erasures.
We limit our analysis to the case where the signal covariance
matrix $\R_{xx}$ is of the form $\R_{xx}=\sigma_x^2\I$. The case
of a general $\R_{xx}$ is more involved and is outside the scope
of this paper.

We prove that maximum robustness against one subspace erasure is
achieved when the fusion frame is tight and all subspaces have
equal dimensions where the optimal dimension depends on SNR. We
also prove that a tight fusion frame consisting of equi-dimensional
subspaces with equal pairwise \textit{chordal distances} is
maximally robust with respect to two and more subspace erasures.
We call such fusion frames \textit{equi-distance tight fusion
frames}. We prove that the pairwise chordal distances between the subspaces
in equi-distance tight fusion frames meet the so-called
\textit{simplex bound}, and thereby establish an intriguing
connection between the construction of such fusion frames and
\textit{optimal Grassmannian packings} (cf. the excellent survey
by Conway et al. \cite{CHS96}). This connection shows that optimal
Grassmannian packings are fundamental for signal processing
applications where low-dimensional projections are used for robust
dimension reduction.

The paper is organized as follows. In Section \ref{sec:LMMSE}, we
derive the MSE in LMMSE estimation of a random vector from its
noisy fusion frame projections. In Section \ref{sec:main}, we
analyze the robustness of LMMSE estimators to erasures of fusion
frame subspaces and derive conditions for the construction of
maximally robust fusion frames. Section
\ref{sec:connectionpackings} establishes a connection between
equi-distance tight fusion frames and optimal Grassmannian packings.
In Section \ref{sec:constructionoptimalff}, we give several
examples for the construction of equi-distance tight fusion
frames. Conclusions are drawn in Section \ref{sec:conclusions}.

\section{LMMSE Estimation from Fusion Frame Measurements}
\label{sec:LMMSE}

Let $\x\in \mathbb{R}^{M}$ be a zero-mean random vector with
covariance matrix $E[\x\x^T]=\R_{xx}$, which we wish to estimate
from $N$ noisy measurement vectors in a fusion frame
$\{\mathcal{W}_i\}_{i=1}^{N}$, with bounds $A\le B$, in the
presence of noise. In other words, we wish to estimate $\x$ from
its noisy projections onto the subspaces
$\{\mathcal{W}_i\}_{i=1}^{N}$. We take the dimension of the $i$th
fusion frame subspace $\mathcal{W}_i$, $i=1,2,\ldots,N$ to be
$m_i$.

Let $\z_i\in \mathbb{R}^M$, $i=1,\ldots,N$ be the measurement
vectors corresponding to $\{\mathcal{W}_i\}_{i=1}^{N}$. The
measurement model for the $i$th subspace $\mathcal{W}_i$,
$i=1,\ldots,N$ is of the form
\[
\z_i=\P_{i}\x+\n_i,
\]
where $\P_{i}\in \mathbb{R}^{M\times M}$ is the orthogonal
projection matrix onto the $m_i$-di\-men\-sio\-nal subspace
$\mathcal{W}_i$, and $\n_i\in \mathbb{R}^{M}$ is the corresponding
noise vector. Assume that the noise vectors for different
subspaces are mutually uncorrelated, and that each noise vector is
white with covariance matrix $E[\n_i\n_i^T]=\sigma_n^2\I$,
$i=1,\ldots,N$. Further, assume that the signal vector $\x$ and
the noise vectors $\n_i$, $i=1,\ldots, N$ are uncorrelated.

Further, define the composite measurement vector $\z\in \mathbb{R}^{NM}$
and the composite projection matrix $\P\in \mathbb{R}^{NM\times
M}$ as $\z=(\z_1^T \ \ \z_2^T \ \ \cdots \ \ \z_N^T)^T$ and
$\P=(\P_1^T \ \ \P_2^T \ \ \cdots \ \ \P_N^T)^T$. Then, the
composite covariance matrix between $\x$ and $\z$ can be written
as
\[
E\left[\begin{pmatrix}\x \\ \z \end{pmatrix}\begin{pmatrix} \x^T &
\z^T
\end{pmatrix}\right]=\begin{pmatrix} \R_{xx} \ & \ \R_{xz}\\ \R_{zx} \
& \ \R_{zz}\end{pmatrix}\in \mathbb{R}^{(N+1)M\times (N+1)M},
\]
where
\[
\R_{xz}=E[\x\z^T]=\R_{xx}\P^T=\R_{xx}\begin{pmatrix}\P_1^T& \cdots
& \P_N^T
\end{pmatrix}
\]
is the $M\times NM$ cross-covariance matrix between $\x$ and $\z$,
$\R_{zx}=\R_{xz}^T$, and
\begin{equation}\label{eq:Rzz}
\R_{zz}=E[\z\z^T]=\P\R_{xx}\P^T+\sigma_n^2\I=\begin{pmatrix}\P_1 \\ \vdots \\
\P_N\end{pmatrix}\R_{xx}
\begin{pmatrix} \P_1^T & \cdots &
\P_N^T\end{pmatrix}+\sigma_n^2\I
\end{equation}
is the $NM\times NM$ composite measurement covariance matrix.

We wish to minimize the MSE in linearly estimating $\x$ from $\z$.
The linear MSE minimizer is known to be the Wiener filter or the
LMMSE filter $\F=\R_{xz}\R_{zz}^{-1}$, which estimates $\x$ by
$\hat{\x}=\F\z$, e.g., see \cite{Scharf-91}. The error covariance
matrix $\R_{ee}$ in this estimation is given by
\begin{align*}
\R_{ee}=E[\e\e^T]&=E\left[(\x-\F\z)(\x-\F\z)^T\right]\nonumber\\
&=\R_{xx}-\R_{xz}\R_{zz}^{-1}\R_{zx}\nonumber\\
&=\R_{xx}-\R_{xx}\P^T(\P\R_{xx}\P^T+\sigma_n^2\I)^{-1}\P\R_{xx}\nonumber\\
&=\Big(\R_{xx}^{-1}+\frac{1}{\sigma_n^2} \P^T\P\Big)^{-1},
\end{align*}
where the last equality follows from the matrix inversion lemma
(Sherman-Morrison-Woodbury formula) \cite[p.50]{Golub-96}.

The MSE is obtained by taking the trace of $\R_{ee}$. Let
$\phi_i$, $i=1,2,\ldots,M$ be the $i$th eigenvalue of
$\R_{xx}^{-1}+(1/\sigma_n^2)\P^T\P$ and assume $\phi_1\ge
\phi_2\ge \cdots \ge \phi_M>0$. Then, the MSE is
\[
MSE=\tr[\R_{ee}]=\sum\limits_{i=1}^M\frac{1}{\phi_i}.
\]
Let $0<\lambda_1\le \lambda_2\le \cdots \le\lambda_M$ be the
eigenvalues of $\R_{xx}$. Then, from \eqref{eq:ff2}, it follows
that
\[
\frac{1}{\lambda_i}+\frac{A}{\sigma_n^2}\le \phi_i \le
\frac{1}{\lambda_i}+\frac{B}{\sigma_n^2}
\]
or alternatively
\[
\frac{1}{\frac{1}{\lambda_i}+\frac{B}{\sigma_n^2}}\le
\frac{1}{\phi_i} \le
\frac{1}{\frac{1}{\lambda_{i}}+\frac{A}{\sigma_n^2}}.
\]
Therefore, we have the following lower and upper bounds for the
MSE:
\[
\sum\limits_{i=1}^M\frac{1}{\frac{1}{\lambda_i}+\frac{B}{\sigma_n^2}}\le
\left(MSE=\sum\limits_{i=1}^M\frac{1}{\phi_i}\right) \le
\sum\limits_{i=1}^{M}\frac{1}{\frac{1}{\lambda_{i}}+\frac{A}{\sigma_n^2}}.
\]
The lower bound is achieved when the fusion frame is tight. That
is, when $A=B$ and
\begin{equation}\label{eq:tff}
\sum\limits_{\ell=1}^N\P_\ell=A\I.
\end{equation}
Taking the trace from both sides of \eqref{eq:tff} yields the
bound $A$ as
\begin{equation}\label{eq:B}
A=\frac{1}{M} \sum\limits_{\ell=1}^N m_{\ell}.
\end{equation}
Thus, the MSE is given by
\begin{equation}\label{eq:MSEg}
MSE=\sum\limits_{i=1}^M\frac{\sigma_n^2\lambda_i}{\sigma_n^2+\frac{\lambda_i}{M}\sum\limits_{\ell=1}^N
m_\ell}.
\end{equation}

\begin{remark}
{\rm When $\R_{xx}=\sigma_x^2\I$, the MSE expression
in \eqref{eq:MSEg} reduces to
\begin{equation}\label{eq:ReeI}
MSE=\frac{M\sigma_n^2\sigma_x^2}{\sigma_n^2+\frac{\sigma_x^2}{M}\sum\limits_{\ell=1}^N
m_\ell}.
\end{equation}}
\end{remark}

\section{Robustness to Subspace Erasures}
\label{sec:main}

We now consider the case where subspace erasures occur, that is
when measurement vectors from one or more subspaces are lost or
discarded. We wish to determine the MSE when the LMMSE filter
$\F$, which is calculated based on the full composite covariance
matrix in \eqref{eq:Rzz}, is applied to the composite measurement
vector with erasures. We do not wish to recalculate the LMMSE
filter every time an erasure occurs. Recalculating the LMMSE
filter requires calculating the inverse of the composite
covariance matrix of the remaining measurement vectors, which is
intractable from a computational standpoint.

In this section, we show how the subspaces in the fusion frame
$\{\mathcal{W}_i\}_{i=1}^N$ must be selected so that the MSE is
minimized under subspace erasures. In our analysis we assume that
$\{\mathcal{W}_i\}_{i=1}^N$ is tight with bound $A$ given by
\eqref{eq:B}. For the sake of simplicity, we limit our analysis to
the case where the signal covariance matrix is
$\R_{xx}=\sigma_x^2\I$. The case where $\R_{xx}$ is a general
positive definite matrix is more involved and is outside the scope
of this paper.

Let $\mathbb{S}\subset \{1,2,\ldots,N\}$ be the set of indices
corresponding to the erased subspaces. Then, the composite
measurement vector with erasures $\tilde{\z}\in \mathbb{R}^{NM}$
may be expressed as
\[
\tilde{\z}=(\I-\E)\z,
\]
where $\E$ is an $NM\times NM$ block-diagonal erasure matrix whose
$i$th $M\times M$ diagonal block is a zero matrix, if $i\notin
\mathbb{S}$, or an identity matrix, if $i\in \mathbb{S}$. In other
words, in $\tilde{\z}$ the measurement vectors associated with the
erased subspaces are set to zero.

The estimate of $\x$ is given by $\tilde{\x}=\F\tilde{\z}$, where
$\F=\R_{xz}\R_{zz}^{-1}$ is the (no-erasure) LMMSE filter. The
error covariance matrix $\widetilde{\R}_{ee}$ for this estimate is
given by
\begin{align*}
\widetilde{\R}_{ee}&=E\left[(\x-\tilde{\x})(\x-\tilde{\x})^T\right]\\
&=E\left[(\x-\F(\I-\E)\z)(\x-\F(\I-\E)\z)^T\right]\nonumber\\
& =\R_{xx}-
\R_{xz}\R_{zz}^{-1}(\I-\E)\R_{zx}-\R_{xz}(\I-\E)^T\R_{zz}^{-1}\R_{zx}\\
&\ \ \ \
+\R_{xz}\R_{zz}^{-1}(\I-\E)\R_{zz}(\I-\E)^T\R_{zz}^{-1}\R_{zx}.
\end{align*}
We can rewrite $\widetilde{\R}_{ee}$ as
\[
\widetilde{\R}_{ee}=\R_{ee}+\overline{\R}_{ee},
\]
where $\R_{ee}=\R_{xx}-\R_{xz}\R_{zz}^{-1}\R_{zx}$ is the
no-erasure error covariance matrix, and
\[
\overline{\R}_{ee}=\R_{xz}\R_{zz}^{-1}\E\R_{zz}\E^T\R_{zz}^{-1}\R_{zx}
\]
is the extra covariance matrix due to erasures. The MSE is given
by
\[
MSE=\tr[\R_{ee}]=MSE_0+\overline{MSE},
\]
where $MSE_0=\tr[\R_{ee}]$ is the no-erasure MSE in
\eqref{eq:ReeI} and
\begin{eqnarray*}
\lefteqn{\overline{MSE} = \tr[\overline{\R}_{ee}]}\\[1ex]
&=&\tr[\R_{xz}\R_{zz}^{-1}\E\R_{zz}\E^T\R_{zz}^{-1}\R_{zx}]\\[1ex]
&=&\tr[\sigma_x^4
\P^T(\sigma_x^2\P\P^T+\sigma_n^2\I)^{-1}\E(\sigma_x^2\P\P^T+\sigma_n^2\I)\E^T(\sigma_x^2\P\P^T+\sigma_n^2\I)^{-1}\P]
\end{eqnarray*}
is the extra MSE due to erasures.

From the matrix inversion lemma \cite[p.50]{Golub-96}, we have
\begin{align}\label{eq:inv1}
(\sigma_x^2\P\P^T+\sigma_n^2\I)^{-1}
&=\frac{1}{\sigma_n^2}\I-\frac{1}{\sigma_n^4}\P(\frac{1}{\sigma_n^2}\P^T\P+\frac{1}{\sigma_x^2}\I)^{-1}\P^T\nonumber\\
&=\frac{1}{\sigma_n^2}\I-\frac{1}{\sigma_n^2}\frac{\sigma_x^2}{A\sigma_x^2+\sigma_n^2}\P\P^T,
\end{align}
where the second equality follows by using
$\P^T\P=\sum_{i=1}^{N}\P_i=A\I$.

Using \eqref{eq:inv1}, we can simplify the expression for
$\overline{MSE}$ to
\begin{align}\label{eq:XMSE1}
\overline{MSE}=\tr[\overline{\R}_{ee}]&=\alpha^2\tr[\P^T\E(\sigma_x^2\P\P^T+\sigma_n^2\I)\E\P]\nonumber\\
& =\alpha^2 \tr\left[\sigma_x^2\left(\sum\limits_{i\in
\mathbb{S}}\P_i\right)^2+\sigma_n^2\left(\sum\limits_{i\in
\mathbb{S}}\P_i\right)\right],
\end{align}
where $\alpha=\sigma_x^2/(\sigma_x^2A+\sigma_n^2)$. The last
equality in \eqref{eq:XMSE1} follows by considering the action of
the erasure matrix $\E$.

We now show how the subspaces in the fusion frame
$\{\mathcal{W}_i\}_{i=1}^N$ must be constructed so that the total
MSE is minimized for a given number of erasures. We consider three
scenarios: one subspace erasure, two subspace erasures, and more
than two subspace erasures.

\subsection{One Subspace Erasure}

If only one of the subspaces, say the $i$th subspace, is erased,
then $MSE$ is given by
\begin{align}
MSE&=MSE_0+\overline{MSE}\nonumber\\
&=MSE_0+\tr[\alpha^2(\sigma_x^2+\sigma_n^2)\P_i]\nonumber\\
&=\frac{M\sigma_x^2\sigma_n^2}{\sigma_n^2 +\frac{\sigma_x^2}{M}\sum\limits_{\ell=1}^N
m_\ell}
+\frac{\sigma_x^4(\sigma_x^2+\sigma_n^2)}{\Big(\sigma_n^2 +\frac{\sigma_x^2}{M}\sum\limits_{\ell=1}^N
m_\ell\Big)^2}\,m_i,\label{eq:MSE25}
\end{align}
where $m_i=\tr[\P_i]$ is the dimension of the $i$th subspace
$\mathcal{W}_i$.

The erasure can occur for any of the subspaces. Thus, we have to
choose $m_i=m$ for all $i=1,\ldots, N$ so that any one-erasure
results in the same amount of performance degradation. This
reduces the MSE expression \eqref{eq:MSE25} to
\[
MSE=\frac{M\sigma_x^2\sigma_n^2}{\left(N
m\sigma_x^2/M+\sigma_n^2\right)}+\frac{\sigma_x^4(\sigma_x^2+\sigma_n^2)m}{\left(N
m\sigma_x^2/M+\sigma_n^2\right)^2}.
\]

As a function of $m$, $MSE=MSE(m)$ has a maximum at
$m=\tilde{m}$, where
\[
\tilde{m}=\frac{M}{N}\frac{(N-1)\sigma_n^4-\sigma_x^2\sigma_n^2}{((N+1)\sigma_n^2+\sigma_x^2)(1-2\sigma_x^2)}.
\]
The MSE is monotonically increasing for $m<\tilde{m}$ and
monotonically decreasing for $m>\tilde{m}$. The smallest value $m$
can take under the constraint that the set of $m$-dimensional
subspace $\{\mathcal{W}_i \}_{i=1}^M$ remains a tight fusion frame
is $m_{min}=\lceil M/N\rceil$, where $\lceil\cdot\rceil$
denotes integer ceiling. We take the largest value $m$ can take to
be $m_{max}\le M$. The maximum allowable dimension $m_{max}$ is
determined by practical considerations. In the distributed sensing
problem it is the maximum number of sensors we can deploy in a
cluster. In the parallel processing problem it is determined by
the maximum computational load that the local processors can
handle, and in the packet encoding problem it corresponds to the
maximum amount of new information (minimum amount of redundancy)
we can include in a packet, while achieving an error correction
goal.

We have the following theorem.

\begin{theorem} The MSE due to the erasure of one subspace is
minimized when all subspaces in $\{\mathcal{W}_i\}_{i=1}^{N}$ have
equal dimension $m=m^\ast$, where
\begin{equation*}
m^\ast=\left\{\begin{array}{*{20}l} m_{min}, & \text{if
$m_{max}\le \tilde m$ or}\\
& \text{if $m_{min}\le
\tilde m\le m_{max}$ and $MSE(m_{min})\le MSE(m_{max})$},\\
m_{max}, & \text{otherwise.}\end{array}\right.
\end{equation*}
\end{theorem}

\subsection{Two Subspace Erasures}

When two subspaces, say the $i$th subspace and the $j$th subspace,
are erased or discarded, the total MSE is given by
\begin{align*}
MSE&=MSE_0+\overline{MSE}\\
&=MSE_0+\alpha^2\tr[\sigma_x^2(\P_i+\P_j)^2+\sigma_n^2(\P_i+\P_j)].
\end{align*}
We take the dimension of all subspaces to be equal to a given $m$
in order to fix the performance against one subspace erasures.
This fixes $MSE_0$ and reduces the minimization of $MSE$ to
minimizing the extra MSE, which is given by
\begin{align*}
\overline{MSE}=2\alpha^2(\sigma_x^2+\sigma_n^2)m+2\alpha^2\sigma_x^2\tr[\P_i\P_j].
\end{align*}

To minimize $\overline{MSE}$ we have to choose $\mathcal{W}_i$ and
$\mathcal{W}_j$, so that $\tr[\P_i\P_j]$ is minimized. Since
$\P_i$ and $\P_j$ are orthogonal projection matrices onto
$\mathcal{W}_i$ and $\mathcal{W}_j$, the eigenvalues of $\P_i\P_j$
are squares of the cosines of the {\em principal angles}
$\theta_\ell(i,j)$, $\ell=1,\ldots,M$ between $\mathcal{W}_i$ and
$\mathcal{W}_j$. Therefore,
\begin{equation} \label{eq:chordal}
\tr[\P_i\P_j]=\sum\limits_{\ell=1}^{M}
\cos^2\theta_\ell(i,j)=M-d_c^2(i,j),
\end{equation}
where
\[
d_c(i,j)=\left(\sum\limits_{\ell=1}^M
\sin^2\theta_\ell(i,j)\right)^{1/2}
\]
is known as the {\em chordal distance} \cite{CHS96} between
$\mathcal{W}_i$ and $\mathcal{W}_j$.

Thus, we need to maximize the chordal distance $d_c(i,j)$. Since
this has to be done for any two subspace erasures, i.e., for any
 pair $(i,j)$, $i\neq j$, we have to construct the subspaces
$\{\mathcal{W}_i\}_{i=1}^N$ so that any such pair has maximum
chordal distance.

In Section \ref{sec:connectionpackings}, we will prove that the subspaces
in a fusion frame consisting of equi-dimensional and equi-distance
(equi-chordal distance) subspaces have maximal chordal distance if
and only if the fusion frame is tight. We call such a fusion frame
an \textit{equi-distance tight fusion frame} and the subspaces
corresponding to it \textit{maximal equi-distance subspaces}. We
note that maximal equi-distance does not mean that the principal
angles between any pair of subspaces must be equal. Therefore,
this is a more relaxed requirement than the equi-isoclinic
condition in \cite{Bod07}.

We have the following theorem.

\begin{theorem}
The MSE due to two subspace erasures is
minimized when the $m$-dimensional subspaces in the tight fusion
frame $\{\mathcal{W}_i\}_{i=1}^N$ are maximal equi-distance
subspaces.
\end{theorem}

We defer the construction of maximal equi-distance subspaces to
Section \ref{sec:connectionpackings}, where we explain the
connection between this construction and the problem of optimal
packing of $N$ planes in a Grassmannian space \cite{CS93, CHS96}.

\subsection{More Than Two Subspace Erasures}

We now consider the case where more than two subspaces are erased
or discarded. Let the subspaces $\{\mathcal{W}_i\}_{i=1}^N$
have equal dimension $m\ge \lceil M/N \rceil$ and equal pairwise
chordal distance $d_c$, so as to fix the performance against one- and
two-erasures. Then, $\overline{MSE}$ can be written as
\begin{eqnarray*}
\overline{MSE}
&=&\alpha^2
\tr\left[\sigma_x^2\left(\sum\limits_{i\in
\mathbb{S}}\P_i\right)^2+\sigma_n^2\left(\sum\limits_{i\in
\mathbb{S}}\P_i\right)\right]\nonumber\\
&=&\alpha^2(\sigma_x^2+\sigma_n^2)\sum\limits_{i\in
\mathbb{S}}\tr[\P_i]+\alpha^2\sigma_x^2\sum\limits_{i\in
\mathbb{S}}\sum\limits_{j\in \mathbb{S},j\neq i}\tr[\P_i \P_j]\\
&=&\alpha^2(\sigma_x^2+\sigma_n^2)|\mathbb{S}|m+\alpha^2\sigma_x^2|\mathbb{S}|(|\mathbb{S}|-1)(M-d_c^2).
\end{eqnarray*}
Similar to the two-erasure case, $\overline{MSE}$ is minimized
when $d_c^2$ takes its maximum value. Thus, we have the following
theorem.

\begin{theorem}
The MSE due to $k$ erasures, $3\le k< N$ is
minimized when the $m$-dimensional subspaces in the tight fusion
frame $\{\mathcal{W}_i\}_{i=1}^N$ are maximal equi-distance
subspaces.
\end{theorem}

\section{Connections  between Tight Fusion Frames and Optimal Packings}
\label{sec:connectionpackings}

In this section, we show that tight fusion frames that consist of
equi-di\-men\-sio\-nal and equi-distance subspaces are closely
related to optimal packings of subspaces. We start by reviewing
the classical packing problem for subspaces \cite{CS93, CHS96}.

\textbf{Classical Packing Problem.} For given $m,M,N$, find a set
of $m$-dimensional subspaces $\{\mathcal{W}_i\}_{i=1}^N$ in
$\RR^M$ such that $\min_{i \neq j} d_c(i,j)$ is as large as
possible. In this case we call $\{\mathcal{W}_i\}_{i=1}^N$ an {\em
optimal packing}.

This problem was reformulated by Conway et al. in \cite{CHS96} by
describing $m$-dimensional subspaces in $\RR^M$ as points on a
sphere of radius $\frac12 (M-1)(M+2)$. This usually provides a
lower-dimensional representation than the Pl\"ucker embedding.
This idea was then used to prove the optimality of many new
packings by employing results from sphere packing theory such as
Rankin bounds for spherical codes. In what follows, we briefly
describe the embedding of the Grassmannian manifold $G(m,M)$ of
$m$-dimensional subspaces of $\RR^M$, as it was described in
\cite{CHS96}. The basic idea is to identify an $m$-dimensional
subspace $\mathcal{W}$ with the traceless part of the projection
matrix $\Q_1$ associated with $\mathcal{W}$, i.e., with
$\overline{\Q}_1 = \Q_1 - \frac{m}{M} \I$. This yields an
isometric embedding of $G(m,M)$ into the sphere of radius
$\sqrt{\frac{m(M-m)}{M}}$ in $\RR^{\frac12 (M-1)(M+2)}$, where the
distance measure is the chordal distance between two projections.
The chordal distance $d_c(\Q_1,\Q_2)$ between two projection
matrices $\Q_1$ and $\Q_2$ is given by
$d_c(\Q_1,\Q_2)=\frac{1}{\sqrt{2}}\|\Q_1-\Q_2\|_2$, and is equal
to $\frac{1}{\sqrt{2}}$ times the straight-line distance between
the projection matrices. This is the reason that $d_c(\Q_1,\Q_2)$
is called {\em chordal} distance. Conway et al. \cite{CHS96}
deduced from this particular embedding the following result.

\begin{theorem}{\rm \cite{CHS96}} Each packing of $m$-dimensional
subspaces $\{\mathcal{W}_i\}_{i=1}^N$ in $\RR^M$ satisfies
\[ d_c^2(i,j) \le \frac{m(M-m)}{M} \frac{N}{N-1}, \quad i,j=1,\ldots,N.\]
\end{theorem}

The upper bound is referred to as the {\em simplex bound}. The
above theorem implies that if the pairwise chordal distances
between a set of $m$-dimensional subspaces of $\RR^M$ meet the
simplex bound those subspaces form an optimal packing, as the
minimum of chordal distances cannot grow any further.

We now establish a connection between tight fusion frames and
optimal packings.

\subsection{Equi-Dimensional Subspaces}

Consider a tight fusion frame $\{\mathcal{W}_i\}_{i=1}^N$, with
bound $A$, consisting of $N$ $m$-di\-men\-sio\-nal subspaces that
do not necessarily have equal pairwise chordal distances. Since
$\{\mathcal{W}_i\}_{i=1}^N$ is tight, we have
\begin{equation}\label{eq:tightness}
A \I = \sum_{i=1}^N \P_{i}.
\end{equation}
On the one hand, we can apply the trace and employ the fact that
$\tr[\P_{i}]=m$ for each $i$, to obtain
\begin{equation}\label{eq:estimateforA1}
AM = Nm.
\end{equation}
On the other hand, we can multiply \eqref{eq:tightness} from left
by $\P_j$ to get
\[ \left(A-1\right)\P_j = \sum_{i=1, i \neq j}^N \P_{j} \P_{i}, \quad j =1,\ldots,N.\]
We can then take the trace, employ the fact that $\tr[\P_{j}]=m$
for each $j$, and use \eqref{eq:chordal}, to obtain
\begin{equation}\label{eq:estimateforA2}
\left(A-1\right)m = \sum_{i=1, i \neq j}^N \tr[\P_j \P_i]
= (N-1)m - \sum_{i=1, i \neq j}^N d_c^2(i,j).
\end{equation}
Equations \eqref{eq:estimateforA1} and \eqref{eq:estimateforA2}
together prove the following result concerning the value of the
fusion frame bound.

\begin{proposition}\label{prop:valueofbound} A tight fusion frame
$\{\mathcal{W}_i\}_{i=1}^N$ with bound $A$ and $m$-di\-men\-sio\-nal
subspaces satisfies
\[ A = \frac{Nm}{M} = N - \sum_{i=1, i \neq j}^N \frac{d_c^2(i,j)}{m}, \quad  j=1,\ldots,N.\]
\end{proposition}

\subsection{Equi-Dimensional and Equi-Distance Subspaces}

We now turn our attention to tight fusion frames
$\{\mathcal{W}_i\}_{i=1}^N$ consisting of equi-dimensional and
equi-distance subspaces, where the common dimension is $m$ and the
common chordal distance is $d_c$. From Proposition \ref{prop:valueofbound}, it follows that
\[
\frac{Nm}{M} = N - (N-1) \frac{d_c^2}{m}.
\]
Thus, $d_c^2$ is given by
\begin{equation}\label{eq:simplex}
d_c^2=\frac{m(M-m)}{M} \frac{N}{N-1},
\end{equation}
which shows that $d_c^2$ precisely equals the simplex bound.

Next we will study whether this condition is sufficient. That is, we
wish to know whether a fusion frame consisting of equi-dimensional
subspaces whose pairwise chordal distances are equal to the
simplex bound is necessarily tight.
%-----------
%We remark that this would lead to constructions of tight fusion
%frames by using results on optimal packing. In fact, we will show
%that this result holds true, and employ this fact in Section
%\ref{sec:LMMSE} to construct fusion frames which are optimally
%resilient against noise and erasures.
%-----------

Consider a fusion frame $\{\mathcal{W}_i\}_{i=1}^N$, consisting of
$N$ $m$-dimensional subspaces with pairwise chordal distances
$d_c$ equal to the simplex bound. Let $\pi_1,\ldots,\pi_M$ be the
eigenvalues of $\P^T \P=\sum_{i=1}^N \P_i$. Since
$\{\mathcal{W}_i\}_{i=1}^N$ is a fusion frame for $\mathbb{R}^M$,
we have $\pi_\ell>0$, $\ell=1,2,\ldots,M$, and the sum of $\pi_\ell$'s
is given by
\begin{equation}\label{eq:sumlambda}
\sum_{\ell=1}^M \pi_\ell = \tr[\P^T \P] = \sum_{i=1}^N \tr[\P_i] =
Nm.
\end{equation}
The sum of $\pi_i^2$'s can be written as
\begin{eqnarray*}
\sum_{\ell=1}^M \pi_\ell^2 & = & \tr[\P^T \P \P^T \P]\\
%& = & \sum_{i=1}^N \tr(\P \P_i \P^T)\\
%& = & \sum_{i=1}^N \sum_{j=1}^N \tr(\P_j \P_i \P_j)\\
& = & \sum_{i=1}^N \sum_{j=1}^N \tr[\P_i \P_j]\\
& = & \sum_{i=1}^N \sum_{j=1, j \neq i}^N \tr[\P_i \P_j] + \sum_{i=1}^N  \tr[\P_i]\\
& = & N(N-1)(m-d_c^2)+Nm,
\end{eqnarray*}
where the last equality follows from \eqref{eq:chordal}. Inserting
the value of the simplex bound, we obtain
\begin{equation}\label{eq:sumlambdasquare2}
\sum_{\ell=1}^M \pi_\ell^2  = \frac{m^2N^2}{M}.
\end{equation}
To conclude that \eqref{eq:sumlambda} together with
\eqref{eq:sumlambdasquare2} implies tightness of the fusion frame,
we consider the problem of minimizing the function
$\sum_{\ell=1}^M \pi_\ell^2$ under the constraint that
$\pi_1,\ldots,\pi_M$ is a sequence of nonnegative values which sum
up to $\sum_{\ell=1}^M \pi_\ell= Nm$. Using the method of Lagrange
multipliers, we see that the minimum is achieved when all
$\pi_\ell$'s are equal to $Nm/M$. This implies that
\eqref{eq:sumlambda} and \eqref{eq:sumlambdasquare2} can be
simultaneously satisfied only when
\[\pi_1=\cdots=\pi_M=\frac{Nm}{M}.\]
From this relation, it follows that $\{\mathcal{W}_i\}_{i=1}^N$ is
a tight fusion frame. Therefore, we have the following theorem.

\begin{theorem} \label{theo:tightsimplex}
Let $\{\mathcal{W}_i \}_{i=1}^N$ be a fusion
frame of $m$-dimensional subspaces with equal pairwise chordal
distances $d_c$. Then the fusion frame is tight if and only if
$d_c$ equals the simplex bound.
\end{theorem}

An immediate consequence of Theorem \ref{theo:tightsimplex} is as follows.

\begin{corollary}
Equi-distance tight fusion frames are
optimal Grassmannian packings.
\end{corollary}

\section{Construction of Equi-Distance Tight Fusion Frames}
\label{sec:constructionoptimalff}

In this section we present a few examples to illustrate the
richness, but also the difficulty of constructing fusion frames
with special properties such as tightness, equi-dimension, and
equi-distance. The optimal packing of $N$ planes in the
Grassmannian space $G(m,M)$ is a difficult mathematical problem,
the solution to which is known only for special values of $N$,$m$,
and $M$. In fact, even optimal packing of lines ($m=1$) or
equivalently constructing equi-angular lines is a deep
mathematical problem. The reader is referred to \cite{SH03} for a
review of problems which are equivalent to the construction of
equi-angular lines. For the construction of optimal packings with
higher-dimensional subspaces we refer the reader to
\cite{CHS96,CS93,CHRSS99}. We would also like to draw the reader's
attention to N. J. Sloane's webpage \cite{webSloane}, which
includes many examples of Grassmannian packings.

\begin{example}\label{exa:Calderbank}
{\rm As our first example for construction of
equi-distance tight fusion frames, we use a result obtained by
Calderbank et al. \cite{CHRSS99} for construction of optimal
packings. The procedure is as follows. Choose $p$ to be a prime
which is either 3 or congruent to $-1$ modulo 8. Then there exists
an explicit construction which produces a tight fusion frame
$\{\mathcal{W}_i \}_{i=1}^{p(p+1)/2}$ in $\RR^p$ with
\[ m_i = \frac{p-1}{2} \quad \mbox{and} \quad d_c^2= \frac{(p+1)^2}{4(p+2)}
\quad \mbox{for all } i,j=1,\ldots,\frac{p(p+1)}{2},\] where $m_i$
denotes the dimension of the $i$th subspace. From Proposition \ref{prop:valueofbound} it
follows that the bound of this fusion frame equals
\[ A = \frac{p^2-1}{4}.\]

As a particular example of this construction we briefly outline
the equi-distance tight fusion frame we obtain for $p=7$. For this, let $Q =
\{q_i\}_{i=1}^3 = \{1,2,4\}$ denote the nonzero quadratic residues
modulo $7$, and $R=\{3,5,6\}$ the nonresidues. Further, let $\H$ be
a $4 \times 4$ Hadamard matrix, e.g.,
\[\H = \begin{pmatrix}
     1 & 1 & 1 & 1\\
     -1 & 1 & -1 & 1\\
     -1 & -1 & 1 & 1\\
     1 & -1 & -1 & 1
 \end{pmatrix}.\]
Finally, we denote the coordinate vectors in $\mathbb{R}^7$ by
$\e_i$, $0 \le i \le 6$, and set $C = \sqrt{2}$ and $k=3$. Then we
define 4 three-dimensional planes $\L_j$, $1 \le j \le 4$ to be
spanned by the vectors
\[ \e_{q_i} + C\H_{ij}\e_{kq_i},\quad 1 \le i \le 3.\]
For each $\L_j$, we obtain $6$ further planes by applying the
cyclic permutation of coordinates $\e_i \mapsto \e_{(i+1) \mod 7}$.
This yields $28$ three-dimensional planes in $\mathbb{R}^7$, which
form a tight fusion frame with bound $12$. Moreover, the chordal
distance between each pair of them equals $d_c^2=\frac{16}{9}$.

This construction is based on employing properties of special
groups, in this case the Clifford group. We remark that this is
closely related with the construction of error-correcting codes.
}
\end{example}

\begin{example} \label{exa:Eisenstein1}
{\rm  This example considers the construction of an
equi-dis\-tance tight fusion frame for a dimension not covered by
Example \ref{exa:Calderbank} by employing the theory of Eisenstein
integers. More precisely, the subspaces will be generated by the
minimal elements of a special lattice. For this, we let
$\mathcal{E} = \{ a + \omega b : a,b \in \mathbb{Z}\}$ denote the
Eisenstein integers, where $\omega = \frac{-1+i\sqrt{3}}{2}$ is a
complex root of unity. The three-dimensional complex lattice
$E_6^*$ over $\mathcal{E}$ is then defined by its generator matrix
\[ \begin{pmatrix}
     \sqrt{-3} & 0 & 0\\
     1 & -1 & 0\\
     1 & 0 & -1
 \end{pmatrix}.\]
It can be shown that the minimal norm of a non-zero element in
$E_6^*$ is $\frac43$. Out of the set of minimal elements, we now
select the following nine:
\[ (1,-1,0), \; (1,0,-1), \; (0,1,-1), \; (\omega,-1,0),\]
\[(0,\omega,-1), \; (-1,0,\omega), \; (\omega,0,-1),
\; (-1,\omega,0), \; (0,-1,\omega).\] Multiplied by the $6$th
roots of unity, this yield $9$ planes in $\mathbb{C}^3$. Using the
canonical mapping of $\mathbb{C}^3$ onto $\mathbb{R}^6$, e.g.,
$(\omega,-1,0) \mapsto (-\frac12, \frac{\sqrt{3}}{2},-1,0,0,0)$,
we obtain $9$ two-dimensional planes in $\mathbb{R}^6$.

In this example all principle angles between each pair of planes
are in fact equal to $\frac{\pi}{3}$. In particular, the chordal
distance is $d_c^2 = \frac32$, which can easily be seen to satisfy
the simplex bound (cf. \eqref{eq:simplex}). By Theorem \ref{theo:tightsimplex} it now
follows that the fusion frame consisting of these planes is tight,
and Proposition \ref{prop:valueofbound} shows that the frame bound equals $3$.
}
\end{example}

\begin{example}
{\rm  The third example explores the construction of
fusion frames in $\mathbb{R}^8$ by employing a similar strategy as
in Example \ref{exa:Eisenstein1}. However, with this example we wish to illustrate the
need to be particularly meticulous when generating a fusion frame
from minimal vectors of a particular lattice. In fact by using a
similar approach, we will generate a tight fusion frame with
equi-dimensional subspaces, but not equi-distance subspaces,
although with a very distinct set of chordal distances.

For our analysis, we choose the lattice
\begin{eqnarray*}
E_8 & = & \left\{(x_1,\ldots,x_8) : \left(x_i \in \mathbb{Z} \;
\forall \, 1 \le i \le 8 \mbox{ or }
x_i \in \mathbb{Z} + \tfrac12 \; \forall \, 1 \le i \le 8\right)\right.\\
& &  \hspace*{2.8cm}\left. \mbox{ and $\sum_{i=1}^8$} x_i \in 2
\mathbb{Z}\right\},
\end{eqnarray*}
which is again a lattice over the Eisenstein integers $\mathcal{E}
= \{ a + \omega b : a,b \in \mathbb{Z}\}$, $\omega =
\frac{-1+i\sqrt{3}}{2}$. Before studying the minimal vectors in
this lattice, we consider the complex root of unity $\omega =
\frac{-1+i\sqrt{3}}{2}$ which was employed in the construction of
$\mathcal{E}$. We first express $\omega$ in quaternions, which
gives $\omega = \frac12(-1+i+j+k)$. Next we define a matrix $\H$
by choosing as row vectors the coefficients of $\omega$,
$i\omega$, $j\omega$, and $k\omega$, i.e.,
\[ \H = \frac12 \begin{pmatrix}
     -1 & 1 & 1 & 1\\
     -1 & -1 & -1 & 1\\
     -1 & 1 & -1 & -1\\
     -1 & -1 & 1 & -1
 \end{pmatrix}.\]
Form this, we build an $8 \times 8$-matrix by setting
\[\mathbf \Omega =\begin{pmatrix}\H & \mathbf 0\\\mathbf 0 & \H\end{pmatrix}.\]
Realizing that this matrix satisfies $\mathbf\Omega^2+\mathbf\Omega+\I=0$, we can
conclude that scaling a vector $\bv \in E_8$ by an Eisenstein
integer $a+\omega b$ can be rewritten as
\[ (a+\omega b)\bv=a\bv+b\mathbf\Omega \bv.\]
Now we are equipped to generate subspaces by minimal vectors,
whose norm can be computed to equal $2$. The lattice $E_8$ has
$240$ minimal vectors, which we assign to planes in the following
way. We first consider the four minimal vectors
\[(1,-1,0,0,0,0,0,0),\; (1,0,-1,0,0,0,0,0),\]
\[ (1,0,0,-1,0,0,0,0),\; (0,1,-1,0,0,0,0,0)\]
and multiply each of them with
\begin{equation}\label{eq:multiplicators}
\I, \; -\I, \; \mathbf\Omega, \; -\mathbf\Omega, \; \I+\mathbf\Omega, \mbox{ and }
-\I-\mathbf\Omega.
\end{equation}
This procedure generates four sets of six minimal vectors, where
each set generates a two-dimensional plane in $\mathbb{R}^8$.
Noticing that this construction only takes all minimal vectors
which are  of the form $(x_1,x_2,x_3,x_4,0,0,0,0)$ into account,
we can clearly use the same idea to group all minimal vectors of
the form $(0,0,0,0,x_5,x_6,x_7,x_8)$. Summarizing, this
construction provides us with $8$ two-dimensional planes in
$\mathbb{R}^8$ which we denote by $\mathcal{W}_1, \ldots,
\mathcal{W}_8$. Next we consider minimal vectors
$(x_1,\ldots,x_8)$, which have one coordinate out of
$x_1,x_2,x_3,x_4$ and one coordinate out of $x_5,x_6,x_7,x_8$
equal to $-1$ or $1$, the others being equal to zero. Again we
multiply these vectors by the factors given in
\eqref{eq:multiplicators}. We can easily see that this procedure
generates another $32$ two-dimensional planes in $\mathbb{R}^8$,
denoted by  $\mathcal{W}_9, \ldots, \mathcal{W}_{40}$.

Although this construction seems similar to the one on Example \ref{exa:Eisenstein1},
we found it surprising to see that in fact
$\{\mathcal{W}_i\}_{i=1}^{40}$ does constitute an equi-dimension
tight fusion frame, however the subspaces are not equi-distance.
The fusion frame bound can be derived from Proposition \ref{prop:valueofbound} and
equals $10$. Most interestingly, the structure of the chordal
distances is rather distinct. In fact, it can be computed that the
chordal distance between each pair is either $d_c^2=2$ -- which
means that they are orthogonal -- or $d_c^2 = \frac43$.
}
\end{example}

\section{Conclusions}
\label{sec:conclusions}

We considered the linear estimation of a random vector from its
noisy projections onto low-dimensional subspaces constituting a
fusion frame. We proved that -- in the presence of white noise -- the
MSE in such an estimation is minimal when the fusion frame is
tight. We analyzed the effect of subspace erasures on the
performance of LMMSE estimators. We proved that maximum robustness
against one subspace erasures is achieved when the fusion frame is
tight and all subspaces have equal dimensions, where the optimal
dimension depends on the SNR. We also proved that equi-distance tight
fusion frames are maximally robust against two and more than two
subspace erasures. In addition we proved that
equi-distance tight fusion frames are in fact optimal Grassmannian
packings, and thereby showed that optimal Grassmannian packings are
fundamental for signal processing applications where
low-dimensional projections are used for robust dimension
reduction. We presented a few examples for the construction of
equi-distance tight fusion frames and illustrated the interesting
and sometimes challenging nature of such constructions.

\section*{Acknowledgment}

The authors would like to thank Stephen Howard for stimulating
discussions about applications of fusion frames.
G. Kutyniok would like to thank Pete Casazza and Christopher Rozell for
interesting general discussions about fusion frames. She
is also indebted to Minh Do and Richard Baraniuk for a
discussion about applications of fusion frames during the 2007 Von
Neuman Symposium in Snow Bird, Utah.  This work was conducted when G.
Kutyniok was a visitor at PACM, Princeton University, and she
would like to thank PACM for its hospitality and support.

\bibliographystyle{elsart-num}
\bibliography{Frames}

\begin{thebibliography}{10}
\expandafter\ifx\csname url\endcsname\relax
  \def\url#1{\texttt{#1}}\fi
\expandafter\ifx\csname urlprefix\endcsname\relax\def\urlprefix{URL }\fi

\bibitem{CK04}
P.~G. Casazza, G.~Kutyniok, Frames of subspaces, in: Wavelets, frames and
  operator theory, Vol. 345 of Contemp. Math., Amer. Math. Soc., Providence,
  RI, 2004, pp. 87--113.

\bibitem{CKL07}
P.~G. Casazza, G.~Kutyniok, S.~Li, Fusion frames and distributed processing,
  preprint.

\bibitem{RGJ06}
C.~J. Rozell, I.~N. Goodman, D.~H. Johnson, Feature-based information
  processing with selective attention, in: Proc. Int. Conf. Acoust., Speech,
  Signal Process. (ICASSP), Vol.~4, 2006, pp. 709--712.

\bibitem{RJ06}
C.~J. Rozell, D.~H. Johnson, Analyzing the robustness of redundant population
  codes in sonsory and feature extraction systems, Neurocomputing 69 (2006)
  1215--1218.

\bibitem{CKLR07}
P.~G. Casazza, G.~Kutyniok, S.~Li, C.~J. Rozell, Modeling sensor networks with
  fusion frames, preprint.

\bibitem{Bod07}
B.~G. Bodmann, Optimal linear transmission by loss-insensitive packet encoding,
  Appl. Comput. Harmon. Anal. 22~(3) (2007) 274--285.

\bibitem{Scharf-91}
L.~L. Scharf, Statistical Signal Processing, Addison-Wesley, MA, 1991.

\bibitem{GKK01}
V.~K. Goyal, J.~Kova{\v{c}}evi{\'c}, J.~A. Kelner, Quantized frame expansions
  with erasures, Appl. Comput. Harmon. Anal. 10~(3) (2001) 203--233.

\bibitem{CK03}
P.~G. Casazza, J.~Kova{\v{c}}evi{\'c}, Equal-norm tight frames with erasures,
  Adv. Comput. Math. 18~(2-4) (2003) 387--430.

\bibitem{HP04}
R.~B. Holmes, V.~I. Paulsen, Optimal frames for erasures, Linear Algebra Appl.
  377 (2004) 31--51.

\bibitem{BP05}
B.~G. Bodmann, V.~I. Paulsen, Frames, graphs and erasures, Linear Algebra Appl.
  404 (2005) 118--146.

\bibitem{SH03}
T.~Strohmer, R.~W. Heath, Jr., Grassmannian frames with applications to coding
  and communication, Appl. Comput. Harmon. Anal. 14~(3) (2003) 257--275.

\bibitem{CK07}
P.~G. Casazza, G.~Kutyniok, Robustness of fusion frames under erasures of
  subspaces and of local frame vectors, preprint.

\bibitem{CHS96}
J.~H. Conway, R.~H. Hardin, N.~J.~A. Sloane, Packing lines, planes, etc.:
  packings in {G}rassmannian spaces, Experiment. Math. 5~(2) (1996) 139--159.

\bibitem{Golub-96}
G.~H. Golub, C.~F. {Van Loan}, Matrix Computations, 3rd Edition, John Hopkins
  Univ. Press, Baltimore, MD, 1996.

\bibitem{CS93}
J.~H. Conway, N.~J.~A. Sloane, Sphere Packings, Lattices and Groups, 2nd
  Edition, Springer, New York, 1993.

\bibitem{CHRSS99}
A.~R. Calderbank, R.~H. Hardin, E.~M. Rains, P.~W. Shor, N.~J.~A. Sloane, A
  group-theoretic framework for the construction of packings in {G}rassmannian
  spaces, J. Algebraic Combin. 9~(2) (1999) 129--140.

\bibitem{webSloane}
www.research.att.com/{$\sim$}njas/grass/index.html.

\end{thebibliography}

\end{document}